\documentclass[11pt,a4,reqno]{amsart}

\usepackage[margin=3.5cm]{geometry}
\usepackage[english]{babel}
\usepackage[utf8]{inputenc}
\usepackage{fancyhdr}
\usepackage{natbib}
\usepackage{amsmath}

\usepackage[font=small, labelfont=small]{caption}

\RequirePackage{amsthm,amsmath,amsfonts,amssymb}
\usepackage{graphicx}
\usepackage{hyperref}
\usepackage{graphicx}
\usepackage{placeins}
\usepackage{caption}
\usepackage{subcaption}
\numberwithin{equation}{section}

\usepackage{xcolor}

\theoremstyle{definition}

\RequirePackage{amsthm,amsmath,amsfonts,amssymb}
\usepackage{graphicx}
\usepackage{hyperref}
\usepackage{lscape}
\usepackage{placeins}
\usepackage{caption}
\usepackage{subcaption}
\numberwithin{equation}{section}

\theoremstyle{plain}
\theoremstyle{remark} 
\numberwithin{equation}{section} 

\usepackage{multirow}

\title{Nonparametric Inference for an Extropy-Based Divergence Measure}

\author{N\lowercase{aresh} G\lowercase{arg$^{a,b}$}, I\lowercase{sha} D\lowercase{ewan$^a$ and} S\lowercase{udheesh} K\lowercase{umar} K\lowercase{attumannil}$^c$\\
    $^a$S\lowercase{tatistical}  S\lowercase{ciences} U\lowercase{nit},
	I\lowercase{ndian} S\lowercase{tatistical} I\lowercase{nstitute}, D\lowercase{elhi}, I\lowercase{ndia},\\
    $^b$P\lowercase{robabilistic}  M\lowercase{achine} L\lowercase{earning} G\lowercase{roup},
	A\lowercase{alto} U\lowercase{niversity}, E\lowercase{spoo}, F\lowercase{inland},\\
    $^c$A\lowercase{pplied} S\lowercase{tatistics} U\lowercase{nit},
	I\lowercase{ndian} S\lowercase{tatistical} I\lowercase{nstitute}, C\lowercase{hennai}, I\lowercase{ndia}.
    }
\thanks{Corresponding author email: skkkattu@isichennai.res.in}
\date{}

\begin{document}
\maketitle
\begin{abstract}
Survival extropy, which quantifies the uncertainty associated with the remaining lifetime distribution, provides an information-theoretic perspective on survival behavior. We consider a divergence measure based on survival extropy and derive its nonparametric estimators based on U-statistics, empirical distribution functions, and kernel density. Further, we construct confidence intervals for the divergence measure using the jackknife empirical likelihood (JEL) method and the normal approximation method with a jackknife pseudo-value-based variance estimator. A comprehensive Monte Carlo simulation study is conducted to compare the performance of the measure with existing divergence measures. Additionally, we evaluate the finite-sample performance of various estimators for the proposed measure.  The findings highlight the effectiveness of the divergence measure and its estimators in practical applications. Finally, we show how the proposed divergence measure is used to detect the small differences between images in image datasets.
\end{abstract}

\section{Introduction}

Divergence measures are fundamental tools in statistics and information theory, quantifying the dissimilarity between two probability distributions.  They are widely employed in diverse fields such as machine learning, signal processing, and bioinformatics for tasks like model comparison, uncertainty quantification, and performance evaluation.  Classical divergence measures, like the Kullback-Leibler divergence (\cite{kullback1951information}), the Hellinger distance (\cite{hellinger1909neue}), and the Jensen-Shannon divergence (\cite{lin1991divergence}), are well-established in literature. However, these measures often encounter difficulties when applied to lifetime data,  which has applications in economics, finance, actuarial, medical and reliability studies.  In this context, survival extropy serves as an alternative measure that captures the uncertainty inherent in the remaining lifetime distribution, offering a complementary perspective to traditional divergence measures (\cite{lad2015extropy}, \cite{jahanshahi2020cumulative}).

Lifetime data analysis focuses on modeling the time until a specific event occurs, such as death, system failure, or  recurrence of disease. Let $X$ be a non-negative absolutely continuous random variable with density function $f$,  distribution function $F$ and    survival function  $\bar{F}(x) = 1 - F(x)$. It is important to compare survival functions  in biomedical research, epidemiology, and actuarial science for evaluating treatment efficacy, identifying risk factors, and facilitating informed decision-making.  
\\

To qauntify the uncertainty in lifetime, \cite{cox2016practical} introduced a divergence measure, \( D_{CC} \), for comparing survival functions:
\begin{equation} \label{eq:1.1}
    D_{CC} = \int_{0}^{\infty} \left| \bar{F}(x)g(x) - \bar{G}(x)f(x) \right| dx.
\end{equation}
The measure $D_{CC} $ can also be expressed as
\begin{equation} \label{eq:1.2}
    D_{CC} = \int_{0}^{\infty}  \bar{F}(x) \bar{G}(x) \left| r_{G}(x) - r_{F}(x) \right| dx,
\end{equation}
where  $r_{F}(x) $ and $r_{G}(x)$ are the failure rates of $X$ and $Y$, respectively. The measure \( D_{CC} \) effectively captures differences between two survival distributions, but as a function of the absolute value function  which introduces non-differentiability. This complicates the study of various properties of its estimators.\\

To address this limitation, we consider a divergence measure which is designed to handle lifetime data  without and with  censoring more effectively. This measure, based on squared integral differences between the survival functions  is:
\begin{equation} \label{eq:1.3}
    D_w = \int_{0}^{\infty} w(x) \left( \bar{F}(x) - \bar{G}(x) \right)^2 dx,
\end{equation}
where \( w(x) \) is the weight function. The divergence measure \( D_w \) provides a smoother and more interpretable alternative to \( D_{CC} \). Unlike \( D_{CC} \), which involves absolute differences, \( D_w \) uses squared differences, thus avoiding the non-differentiability issues that often arise when dealing with absolute values. This choice not only makes the estimation process easier but also makes the measure easier to understand. Close examination reveals that the measure $D_w$ is based on survival extropy (see below for a detailed discussion), which makes it a suitable candidate for analysing lifetime data.

In the divergence measure \( D_w \), the weight function \( w(x) \) plays an important role in controlling the significance of the differences between the survival functions \( \bar{F}(x) \) and \( \bar{G}(x) \) at various time points. By adjusting \( w(x) \), one can give more importance to specific regions of the survival curves. For example, if early survival differences are more important, a weight function that assigns higher values to smaller values of \( x \) can be chosen. Alternatively, if long-term survival is more important, \( w(x) \) can be designed to give more weight to larger values of \( x \). This flexibility allows the divergence measure to focus on the most relevant differences for the analysis, making it adaptable to the context of the data.

Moreover, if \( w(x) \) is chosen to be constant for all \( x \), i.e., \( w(x) = w \) for all \( x \), it acts as a global scaling factor. This simplifies the divergence measure while still keeping its general form. To make \( D_w \) to be scale-invariant, we can normalize it by selecting
\begin{equation} \label{eq:1.4}
w = \frac{1}{\int_0^\infty (\bar{F}(x) + \bar{G}(x)) \, dx},
\end{equation}
which ensures that \( D_w \) remains bounded regardless of the scale of survival times. Alternatively, one could choose
\begin{equation} \label{eq:1.5}
w = \frac{1}{\int_0^\infty \max(\bar{F}(x), \bar{G}(x))^2 \, dx},
\end{equation}
which keeps \( D_w \) scale-invariant and comparable across different datasets. Both (1.4) and (1.5) normalizations ensure that \( D_w \) always lies within the interval \( [0, 1] \), where:
\begin{itemize}
    \item \(
D_w = 0 \quad \quad\;\; \;   \text{if} \quad \bar{F}(x) = \bar{G}(x) \quad \text{for all} \quad x,
\)
\item \(
0< D_w < 1 \quad  \text{if} \quad \bar{F}(x) \ne \bar{G}(x) \quad \text{for all} \quad x\).
\end{itemize}
By using these normalization techniques, \( D_w \) remains a clear, interpretable, and flexible divergence measure that can be applied to different types of data. For the choice \( w(x) = 1/2 \), \cite{saranya2024relative} discussed about this measure and explored some of its theoretical properties.

  Furthermore, as mentioned above \( D_w \) is closely related to extropy, introduced by \cite{lad2015extropy}, which quantifies the uncertainty inherent in the survival functions. This connection provides a strong theoretical foundation for its application in survival analysis.  Recent work has further explored the relationship between extropy and divergence measures (e.g., \cite{qiu2018extropy}, \cite{yang2019bounds}, \cite{kayal2021failure}, \cite{kharazmi2024fisher}, \cite{kharazmi2024fractional} and \cite{saranya2025inaccuracy}).\\

\noindent Expanding \( D_w \) yields the following:
\begin{equation} \label{eq:1.6}
    D_w = \int_{0}^{\infty} w(x) \bar{F}^2(x) dx + \int_{0}^{\infty}w(x) \bar{G}^2(x) dx - 2 \int_{0}^{\infty}w(x) \bar{F}(x)\bar{G}(x) dx.
\end{equation}\\

The first two terms are related to the cumulative residual entropy of \( X \) and \( Y \), respectively (\cite{jahanshahi2020cumulative}), while the third term is related to the measure of inaccuracy of cumulative residual extropy. This formulation provides deeper insight into the properties of \( D_w \) and its relevance to survival analysis. \\

\noindent For simplicity and better clarity in \( D_w \), we choose \( w(x) = 1 \) throughout the paper, and denote it as $D$. After algebraic manipulations, \( D \) can also be rewritten as:

\begin{equation} \label{eq:1.7}
D=\int_{0}^{\infty} 2 x \bar{F}(x) d F(x)+\int_{0}^{\infty} 2 x \bar{G}(x) d G(x)-\int_{0}^{\infty} 2 x \bar{F}(x) d G(x)-\int_{0}^{\infty} 2 x \bar{G}(x) d F(x).
\end{equation}

In particular, the last two terms also appear in the definition of measure \( D_{CC} \) (see \eqref{eq:1.1}). The representation in \eqref{eq:1.7} enables measure \( D \) to be applied to both the complete data and right-censored data.\\

To estimate the proposed measure, we employ non-parametric approaches such as U-statistics, empirical distribution-based estimators, and kernel-based methods. These techniques offer flexible, data-driven estimation. Through a simulation study, we compare the proposed estimators of $D$ with those of existing divergence measures. The results indicate that the proposed measure performs favorably against traditional measures across various scenarios.  We also construct confidence intervals for $D$ using two approaches: the jackknife empirical likelihood (JEL) method and the normal approximation with jackknife pseudo-values. The JEL method provides a more robust, distribution-free alternative, while the normal approximation approach offers computational simplicity. Our simulation results indicate that the JEL-based confidence intervals often outperform those obtained via normal approximation.  
\\

The paper proceeds as follows: Section 2 presents various non-parametric approaches used to derive estimators of $D$. In Section 3, we construct the confidence intervals for $D$ using both the normal approximation and the jackknife empirical likelihood approach. Section 4 discusses simulation studies that compare the estimator of $D$ with kernel-type estimators of existing measures using the mean squared error (MSE) criterion. We also evaluate the performance of the estimators and the confidence intervals for $D$. Section 5 presents real-life applications of the divergence measure to various datasets, including image datasets.  Finally, Section 6 provides concluding remarks.


\section{Estimation  of $D$}

In this section we propose three nonparametric estimators  of \( D \) based on  (i) U-statistics (ii) empirical distribution function  and (iii) kernel type  density estimators based on   two independent random samples \( X_1, X_2, \dots, X_{n_1} \) and \( Y_1, Y_2, \dots, Y_{n_2} \)  drawn from distributions with distribution functions \( F \) and \( G \), respectively.   Each of these approaches offers unique advantages depending on the nature of the data and underlying assumptions. The U-statistics approach ensures asymptotic efficiency, the empirical approach is straightforward and widely used in practice, while the kernel-based estimator provides smooth estimates even in cases where the data is sparsely distributed.

In Subsections 2.1-2.3, we provide a detailed formulation of the different estimators along with their theoretical justifications and computational implementations. 

\subsection{U-statistics-based estimator}
First, we derive an alternative representation of \( D \). Let \( X_1 \) and \( X_2 \) be two independent random variables with distribution function \( F \), and similarly, let \( Y_1 \) and \( Y_2 \) be two independent random variables with distribution function \( G \). Note that, for a non-negative random variable \( X \), its expectation is given by  
\[
E(X) = \int_{0}^{\infty} \bar{F}(x) \, dx.
\]
Noting that \( \bar{F}^{2}(x) \), \( \bar{G}^{2}(x) \), and \( \bar{F}(x) \bar{G}(x) \) correspond to the survival functions of the random variables \( \min(X_1, X_2) \), \( \min(Y_1, Y_2) \), and \( \min(X_1, Y_1) \) respectively, we can express \( D \) as follows:  
\begin{align}
	D &=\int_{0}^{\infty} (\bar{F}(x)-\bar{G}(x))^2 \, dx \nonumber\\
	&= \int_{0}^{\infty} \bar{F}^2(x) \, dx + \int_{0}^{\infty} \bar{G}^2(x) \, dx - 2 \int_{0}^{\infty} \bar{F}(x)\bar{G}(x) \, dx  \nonumber\\
	&= E[\min\{X_1, X_2\}] + E[\min\{Y_1, Y_2\}] - 2E[\min\{X_1, Y_1\}].  \label{eq:2.1}
\end{align}  

Next, we construct a U-statistics-based estimator of \( D \). Define a symmetric kernel function \( h(X_1, X_2, Y_1, Y_2) \) of degree $(2,2)$  as  
\begin{align} \label{eq:2.2}
	h(X_1,X_2,Y_1,Y_2) &= \frac{1}{4} \Big[ 4 \min\{X_1,X_2\} + 4 \min\{Y_1,Y_2\} - 2 \min\{X_1,Y_1\} \nonumber \\
	&\quad - 2 \min\{X_1,Y_2\} - 2 \min\{X_2,Y_1\} - 2\min\{X_2,Y_2\} \Big].
\end{align}  

The corresponding U-statistic estimator of \( D \) is given by  

\begin{equation}  \label{eq:2.3}
\hat{D}_{Ustat}= \frac{1}{\binom{n_1}{2}\binom{n_2}{2}}	\sum_{1\leq i_1<i_2\leq n_1} \sum_{1\leq j_1<j_2\leq n_2} h(X_{i_1},X_{i_2},Y_{j_1},Y_{j_2}).
\end{equation}  

By definition, estimator $\hat{D}_{Ustat}$ is unbiased for $D$ and it also consistent for \( D \) (see \cite{lehmann1951consistency}). 

Next, we find the asymptotic distribution of $\hat{D}_{Ustat}$. For this purpose, we define the following functions:
\begin{align*}
    g_{10}(x) & = E \left[h(x, X_{2}, Y_{1}, Y_2)\right] - \theta, \quad \sigma_{10}^{2} = \operatorname{var}\left(g_{10}(X_{1})\right), \\
    g_{01}(y) & = E \left[h(X_{1}, X_{2}, y, Y_2)\right] - \theta, \quad \sigma_{01}^{2} = \operatorname{var}\left(g_{01}(Y_{1})\right).
\end{align*}

The following theorem establishes the asymptotic distribution of \( \hat{D}_{Ustat} \).

\vspace*{3mm}

\noindent \textbf{Theorem 1.} Let $n=n_1+n_2$, and suppose $\frac{n_1}{n}\rightarrow p\in (0,1)\;$ as $n\rightarrow \infty$. Then $\sqrt{n} (\hat{D}_{Ustat}-D)$ converges in distribution to Normal r.v. with mean zero and variance $\sigma^2$, where
$$\sigma^2=\frac{4}{p} \sigma^2_{10}+\frac{4}{1-p} \sigma^2_{01},$$
with
\begin{align*}
    \sigma_{10}^2&=Var\left(\int_{0}^{X} z \,dF(z) - \int_{0}^{X} z \, dG(z) + X(\bar{F}(X)-\bar{G}(X)) \right),\\
\text{and}\qquad \sigma_{01}^2&=Var\left(\int_{0}^{Y} z\, dG(z) - \int_{0}^{Y} z \, dF(z) + Y(\bar{G}(Y)-\bar{F}(Y)) \right).
\end{align*}

\vspace{2mm}

\subsection{Empirical distribution function-based estimator}

 Let \( X_{(1)}, X_{(2)}, \dots, X_{(n_1)} \) and \( Y_{(1)}, Y_{(2)}, \dots, Y_{(n_2)} \) be the order statistics corresponding to the random samples \( X_1, X_2, \dots, X_{n_1} \) and \( Y_1, Y_2, \dots, Y_{n_2} \), respectively.  \\

The empirical distribution function \( F_{n_1}(x) \), which serve as estimator of the true CDFs \( F(x) \), is defined as follows:

\begin{equation} \label{eq:2.4}
F_{n_1}(x) =
\begin{cases} 
    0, & \text{if } x \leq X_{(1)} \\
    \frac{j}{n_1}, & \text{if } X_{(j)} < x \leq X_{(j+1)}, \quad j=1,2,\dots,n_1-1 \\
    1, & \text{if } x > X_{(n_1)}
\end{cases}.
\end{equation}
Similarly, \( G_{n_2}(x) \) can be defined based on the order statistics \( Y_{(1)}, Y_{(2)}, \dots, Y_{(n_2)} \).



Here, we propose an empirical estimator for \( D \), as defined in equation \eqref{eq:1.7}.  Let \( S_{(j)} \) and \( R_{(j)} \) represent the ranks of \( X_{(j)} \) and \( Y_{(j)} \), respectively, in the combined increasing arrangement of all \( X \)'s and \( Y \)'s. Then, the empirical distribution function-based estimator of \( D \) is given by
\begin{align} \label{eq:2.5}
\hat{D}_{Emp} = \frac{1}{n_1} \sum_{j=1}^{n_1-1} (X_{(j)} + X_{(j+1)}) &
\left( \frac{S_{(j)}}{n_2} - j \left(\frac{1}{n_1} + \frac{1}{n_2} \right) \right) \nonumber\\
+ \frac{1}{n_2} & \sum_{j=1}^{n_2-1} (Y_{(j)} + Y_{(j+1)}) 
\left( \frac{R_{(j)}}{n_1} - j \left(\frac{1}{n_1} + \frac{1}{n_2} \right) \right).
\end{align}

Estimators using a similar approach have also been considered in the literature by \cite{kochar1981new} and \cite{deshpande1972linear}.

 \subsection{Kernel-based estimator}

 In this section, we consider a non-parametric estimator for the divergence measure \( D \) using kernel density estimation. A similar type of estimator was also discussed by \cite{saranya2024relative}; however, here we provide an estimator for a more general scenario, allowing for different sample sizes and bandwidths.  The kernel density estimator of $f(x)$, introduced by \cite{parzen1962estimation}, is given by:
\begin{equation} \label{eq:2.6}
    \hat{f}_{n}(x) = \frac{1}{n b_{n}} \sum_{j=1}^{n} k\left(\frac{x - X_{j}}{b_{n}}\right),
\end{equation}
where $k(\cdot)$ is a kernel function that satisfies the following conditions: $k(x) \geq 0$ for all $x$, $\int k(x) \,dx = 1$, $k(\cdot)$ is symmetric about zero, and $k(\cdot)$ satisfies the Lipschitz condition. The bandwidth parameter $b_n$ is a sequence of positive numbers such that $b_n \to 0$ and $n b_n \to \infty$ as $n \to \infty$.  \\

\noindent Let \(\bar{K}(t) = \int_{t}^{\infty} k(u) \, du.\) Then, using kernel density estimation, the non-parametric estimator of the survival function is given by:
\begin{equation} \label{eq:2.7}
    \hat{\bar{F}}_{n_1}(x) = \frac{1}{n_1} \sum_{j=1}^{n_1} \bar{K} \left( \frac{x - X_{j}}{b_{n_1}} \right).
\end{equation}

Similarly, the survival function estimator for $G$ is denoted as $\hat{\bar{G}}_{n_2}(x)$. Based on these estimates, the non-parametric kernel estimator for $D$ is defined as:
\begin{align} \label{eq:2.8}
    \hat{D}_{Ker} &= \int_{0}^{\infty} \left( \hat{\bar{F}}_{n_1}(x) - \hat{\bar{G}}_{n_2}(x) \right)^2 \, dx \nonumber \\
    &= \int_{0}^{\infty} \left( \left( \frac{1}{n_1} \sum_{j=1}^{n_1} \bar{K} \left( \frac{x - X_{j}}{b_{n_1}} \right) \right) - \left( \frac{1}{n_2} \sum_{j=1}^{n_2} \bar{K} \left( \frac{x - Y_{j}}{b_{n_2}} \right) \right) \right)^2 dx.
\end{align}

\noindent We use the Gaussian kernel function, which is widely used due to its smoothness and favorable theoretical properties. The bandwidths are computed using Silverman’s rule of thumb (\cite{silverman2018}):
\begin{align*}
    b_{n_1} &= 0.9 \min \left\{ \operatorname{sd}(X), \frac{\operatorname{iqr}(X)}{1.34} \right\} n_1^{-0.2}, \\
    b_{n_2} &= 0.9 \min \left\{ \operatorname{sd}(Y), \frac{\operatorname{iqr}(Y)}{1.34} \right\} n_2^{-0.2},
\end{align*}
where $\operatorname{sd}(\cdot)$ represents the standard deviation, and $\operatorname{iqr}(\cdot)$ denotes the interquartile range. These bandwidths control the degree of smoothing in the kernel density estimation process.

 The kernel-based estimator provides a flexible and data-driven approach to estimating $D$. Using kernel density estimation, this method captures distributional differences in a smooth and efficient manner. The performance of the estimators are evaluated through Monte Carlo Simulation.

\section{Confidence interval for $D$}
In the next sections, we present the construction of confidence intervals for the measure $D$. We obtain JEL and a normal-based confidence interval for $D.$
\subsection{Normal approximation with Jackknife pseudo-values for C.I. of $D$}
In this subsection, we construct a confidence interval for $D$ using the normal approximation (Theorem 1) of U-statistics. As a consequence of Theorem 1, the asymptotic distribution of the U-statistic estimator of $D$ is as follows:
\begin{equation} \label{eq:2.9}
\frac{\hat{D}_{Ustat} - D}{\sigma/\sqrt n} \xrightarrow{d} N(0, 1) \quad \text{as } n \to \infty,
\end{equation}
where \( \xrightarrow{d} \) denotes convergence in distribution.\\

Since \(\sigma^2\) is unknown, we need its consistent estimator. We consider the jackknife pseudo-values-based estimator of \( \sigma^2 \), proposed by \cite{arvesen1969jackknifing} and further discussed by \cite{wang2010empirical}, \cite{guangming2013nonparametric},\cite{an2018jackknife}, and \cite{garg2024jackknife}. The estimator is given by
\begin{align} \label{eq:2.10}
    \hat{\sigma}^{2} &= \frac{1}{n_{1}\left(n_{1}-1\right)} \sum_{i=1}^{n_{1}}\left(V_{i, 0}-\bar{V}_{\cdot, 0}\right)^{2} 
    +\frac{1}{n_{2}\left(n_{2}-1\right)} \sum_{j=1}^{n_{2}}\left(V_{0, j}-\bar{V}_{0,\cdot}\right)^{2},
\end{align}
where
\begin{align*}
    V_{i,0} &= n_{1} U_{n_{1}, n_{2}} - \left(n_{1}-1\right) U_{n_{1}-1, n_{2}}^{(-i)},\quad \forall \; i=1,2,\dots,n_1,\\
    V_{0,j} &= n_{2} U_{n_{1}, n_{2}} - \left(n_{2}-1\right) U_{n_{1}, n_{2}-1}^{(-j)},\quad \forall \; j=1,2,\dots,n_2.
\end{align*}
Here, \( U_{n_{1}-1, n_{2}}^{(-i)} \) is the U-statistic obtained after deleting the observation \( X_{i} \), and \( U_{n_{1}, n_{2}-1}^{(-j)} \) is the U-statistic obtained after deleting \( Y_{j} \). The mean values are computed as
\begin{align*}
    \bar{V}_{\cdot, 0} = \frac{1}{n_1}\sum_{i=1}^{n_1}V_{i,0},\quad \text{and} \quad
    \bar{V}_{0, \cdot} = \frac{1}{n_2}\sum_{j=1}^{n_2}V_{0,j}.
\end{align*}

\cite{arvesen1969jackknifing}, along with \cite{wang2010empirical}, \cite{guangming2013nonparametric}, \cite{an2018jackknife}, and \cite{garg2024jackknife}, showed that \( \hat{\sigma}^{2} \) is a consistent estimator of \( \sigma^{2} \).\\

\noindent Therefore, a two-sided \( (1-\alpha) \) confidence interval for \( D \), based on \( \hat{\sigma}^{2} \), is given by
\begin{equation} \label{eq:Normal App}
    \left(\hat{D}_{Ustat}-z_{\alpha / 2} \hat{\sigma}/\sqrt{n}, \hat{D}_{Ustat}+z_{\alpha / 2} \hat{\sigma}/\sqrt{n}\right),
\end{equation}
where \( z_{\alpha} \) is the upper \( \alpha \)-th percentile point of the standard normal distribution.
\subsection{Jackknife empirical likelihood (JEL) based C.I. of $D$}
 The JEL method, introduced by \cite{jing2009jackknife}, extends the traditional Empirical Likelihood (EL) approach to handle nonlinear statistics more effectively. They demonstrated the superior performance of JEL, particularly for one- and two-sample U-statistics, in comparison to the normal approximation method. Given this advantage, we employ the JEL method to construct confidence intervals for the measure $D$.\\

To implement the JEL method, we first merge the two samples ${X_{1}, \ldots, X_{n_{1}}}$ and ${Y_{1}, \ldots, Y_{n_{2}}}$, into a single combined sample by defining a new variable $W_i$ for $i = 1, 2, \dots, n$ as follows:

\begin{equation} \label{eq:2.11}
W_{i}= \begin{cases}X_{i}, & i=1, \ldots, n_{1} \\ Y_{i-n_{1}}, & i=n_{1}+1, \ldots, n .\end{cases}
\end{equation}
where $n=n_1+n_2$. Define $T_n\left(W_{1}, W_{2},\ldots, W_{n}\right) =\hat{D}_{Ustat},$ where $\hat{D}_{Ustat}$ is a $U$-statistic of degree $(2, 2)$ with a symmetric kernel $h$ (defined in \eqref{eq:2.3}).\\

The jackknife pseudo-values $V_i$ for a combined single sample $(W_1,W_2,\dots,W_n)$ are defined as:
\begin{equation} \label{eq:2.12}
    V_{i}=n T_{n}-(n-1)\, T_{n-1}^{(-i)}, \quad i=1,2,\dots,n,
\end{equation}
where $T_{n-1}^{(-i)}$ represents the value of the U-statistic $\hat{D}_{Ustat}$ computed using the $(n-1)$ observations, excluding $W_i$. Also, it is straightforward to show that the U-statistic can be expressed in terms of the pseudo-values as:
\begin{equation} \label{eq:2.12}
    \hat{D}_{Ustat}=\frac{1}{n} \sum_{i=1}^{n} V_{i},
\end{equation}
where the $V_{i}$'s are asymptotically independent random variables.

Now, to construct confidence intervals for the parametric function $D$ using the empirical likelihood (EL) method, we maximise the empirical likelihood function under constraints based on the jackknife pseudo-values. Define a probability vector $\mathbf{p} = (p_1, p_2, \dots, p_n)$, where each pseudo-value $V_i$ is assigned a probability $p_i$ and let $\theta = D$. The empirical likelihood function evaluated at $\theta$ is given by:
\begin{equation} \label{eq:2.14}
    L(\theta) = \max \left\{ \prod_{i=1}^{n} p_{i} \; \Bigg| \; \sum_{i=1}^{n} p_{i} = 1, \sum_{i=1}^{n} p_{i} (V_{i}-E[V_{i}]) = 0 \right\}.
\end{equation}

The expectation $E[V_k]$ is given by:
\begin{equation} \label{eq:2.15}
    E [V_{k}]=\begin{cases}
        \theta\left(\frac{n}{n-2}\right)\left[\left(n_{2}-1\right) \frac{2}{n_{1}}-1\right], & \text{if } k=1, \ldots, n_{1}, \\
        \theta\left(\frac{n}{n-2}\right)\left[\left(n_{1}-1\right) \frac{2}{n_{2}}-1\right], & \text{if } k=n_{1}+1, \ldots, n.
    \end{cases}
\end{equation}

Since $\prod_{i=1}^{n} p_{i}$ attains its maximum $1/n^{n}$ at $p_{i}=1/n$ under the constraint $\sum_{i=1}^{n} p_{i}=1$, the jackknife empirical likelihood ratio at $\theta$ is given by:
\begin{equation} \label{eq:2.16}
    R(\theta) = \frac{L(\theta)}{n^{-n}} = \max \left\{ \prod_{i=1}^{n} \left(n p_{i}\right) \; \Bigg| \; \sum_{i=1}^{n} p_{i}=1, \sum_{i=1}^{n} p_{i} (V_{i}-E[ V_{i}]) = 0 \right\}.
\end{equation}

Applying the method of Lagrange multipliers, when:
\begin{equation*}
    \min _{1 \leq i \leq n} (V_{i}-E [V_{i}]) < 0 < \max _{1 \leq i \leq n} (V_{i}-E [V_{i}]),
\end{equation*}
we obtain the probability weights:
\begin{equation*}
    p_{i}=\frac{1}{n} \frac{1}{1+\lambda(V_{i}-E [V_{i}])},
\end{equation*}
where $\lambda$ satisfies:
\begin{equation*}
    \frac{1}{n} \sum_{i=1}^{n} \frac{V_{i}-E V_{i}}{1+\lambda(V_{i}-E [V_{i}])} = 0.
\end{equation*}

Thus, the non-parametric jackknife empirical log-likelihood ratio is given by:
\begin{equation} \label{eq:2.17}
    \log R(\theta) = -\sum_{i=1}^{n} \log \left(1+\lambda(V_{i}-E [V_{i}])\right).
\end{equation}

To construct confidence intervals or hypothesis tests for $\theta$, we need the asymptotic distribution of $\log R(\theta)$. The following theorem, analogous to Wilks' theorem, provides this result. The proof follows the same lines as Theorem 2 of \cite{jing2009jackknife} and is omitted here.\\
\textbf{Theorem 2.} Assume that $E[ h^{2}\left(X_{1}, X_2, Y_{1}, Y_2\right)]<\infty,\; \sigma_{10}^{2}>0$ and $\sigma_{01}^{2}>$ 0. Also, $0<\liminf n_{1} / n_{2} \leq \lim \sup n_{1} / n_{2}<\infty$. Then, we have $$-2 \log R(\theta) \xrightarrow{d} \chi_{1}^{2}.$$

\vspace{3mm}

\noindent Using Theorem 2, we construct a confidence interval at the level $(1-\alpha)$ for $\theta$, given by
	\begin{equation} \label{eq: CI JEL}
	    CI =
	\{ \theta\;\vert \;-2 \log R(\theta)\leq \chi^2_{
		1,1-\alpha}\},
	\end{equation} 
	where $\chi^2_{1,1-\alpha}$
	is the $(1 -\alpha)^{th}$ percentile point of a $\chi^2$ distribution with one degree of freedom.

\vspace{3mm}

\section{Monte Carlo Simulations}

In this section, we conduct Monte Carlo simulations to evaluate the performance of the estimator for the measure $D$, comparing it with estimators of two well-established measures. Additionally, we assess the performance of various point estimators for $D$ introduced in Section 2. We use the mean squared error (MSE) criterion to compare the performance of various estimators.\\

For the simulation study, we consider two independent exponential distributions with parameters $\lambda_1$ and $\lambda_2$. Independent random samples of sizes $n_1$ and $n_2$ are drawn from these distributions. In Sections 3.1-3.2, The simulation procedure is repeated 2,000 times to ensure robust results and reliable statistical inference.

\subsection{Comparison of MSEs for $D$ and estimators of existing measures }

In this section, we compare the mean squared error (MSE) of the estimator for the measure $D$ with the estimators of $D_{CC}$ (introduced by \cite{cox2016practical}) and the Kullback-Leibler (KL) divergence. The kernel-based estimator for $D$ was derived in Section 2.3. The kernel-based estimators for $D_{CC}$ and the KL divergence are given as follows:

\begin{equation} \label{eq:DCC est}
\hat{D}_{CC} = \int_{0}^{\infty} \left| \hat{\bar{F}}_{n_1}(x) \hat{g}_{n_2}(x) - \hat{\bar{G}}_{n_2}(x) \hat{f}_{n_1}(x) \right| \,dx,
\end{equation}

\begin{equation} \label{eq:KL est}
\hat{KL} = \int_{0}^{\infty} \hat{f}_{n_1}(x) \log\left(\frac{\hat{f}_{n_1}(x)}{\hat{g}_{n_2}(x)}\right) \,dx,
\end{equation}
where 
$\hat{f}_{n_1}(x), \hat{g}_{n_2}(x), \hat{\bar{F}}_{n_1}(x), \hat{\bar{G}}_{n_2}(x), b_{n_1}$ and
$b_{n_2}$ are as defined in Section 2.3.\\

The exact expressions for these three divergence measures in the case of two exponential distributions with rate parameters \( \lambda_1 \) and \( \lambda_2 \) are given by:  

\begin{align*}  
    D = \frac{1}{2\lambda_1} + \frac{1}{2\lambda_2} - \frac{2}{\lambda_1+\lambda_2}, \;\;
    D_{CC} = \frac{|\lambda_1 - \lambda_2|}{\lambda_1 + \lambda_2}, \;\;
    KL = \log\left(\frac{\lambda_1}{\lambda_2}\right) + \left(\frac{\lambda_2}{\lambda_1}\right) - 1.  
\end{align*}

\noindent To compare the performance of three estimators $\hat{D}_{Ker}$, $\hat{D}_{CC}$, and $\hat{KL}$, we calculated the MSEs relative to the true values across various sample sizes and parameter settings. The results are summarized in Table 1. Based on these results, we make the following key observations:

1. The MSEs of the estimators for all three divergence measures consistently decrease as the sample size increases. This aligns with theoretical expectations, as larger sample sizes provide more reliable and stable estimates.
    
    2. The MSE of $\hat{D}_{Ker}$ is consistently less than that of $\hat{D}_{CC}$ and $\hat{D}_{KL}$ across the various choices of $(n_1, n_2)$ and $(\lambda_1, \lambda_2)$ considered. This indicates that $\hat{D}_{Ker}$ provides a more efficient and accurate estimation compared to $\hat{D}_{CC}$ and $\hat{D}_{KL}$.

 Moreover, beyond its empirical advantages in estimation accuracy, $D$ itself offers several theoretical benefits, such as differentiability and smoothness. These properties make it a more preferable choice in various statistical and practical applications, especially where smoothness and analytical tractability are essential.

 Overall, the findings highlight that $D$ not only competes well with existing divergence measures but also provides additional advantages, making it a strong candidate for broader use in statistical inference and machine learning applications.

\subsection{Comparison of estimators for $D$}

To evaluate the performance of different estimators for $D$, we computed estimates using the kernel-based, empirical-based, and U-statistics-based methods. The MSEs of these estimators were calculated across various sample sizes and parameter settings, with the results summarized in Table 2. Based on these results, we make the following key observations:
\\

\FloatBarrier
\begin{table}[]
\caption{Comparison of MSE of estimators of three measures for various sample of sizes $(n_1,n_2)$ drawn from two exponential distributions with parameters ($\lambda_1$,$\lambda_2$).}
\begin{tabular}{ccccc}
\hline
\multicolumn{5}{c}{{($\lambda_1$,$\lambda_2$)=(0.1,0.1)}}                                                                           \\ \hline
Measure / $(n_1,n_2)$ & (10,10)    & (50,40)    & (80,100)   & (200,200)  \\ \hline
$\hat{D}_{Ker}$                 & 1.2084     & 0.0666     & 0.0169     & 0.0034     \\
$\hat{D}_{CC}$            & 0.1021     & 0.0345     & 0.0202     & 0.0107     \\
$\hat{KL}$                    & 0.4330     & 0.0244     & 0.0077     & 0.0029     \\ \hline
\multicolumn{5}{c}{ {($\lambda_1$,$\lambda_2$)=(0.5,1)}}                              
                                                     \\ \hline
Measure / $(n_1,n_2)$ & (10,10)    & (50,40)    & (80,100)   & (200,200)  \\ \hline
$\hat{D}_{Ker}$                & 0.0764     & 0.0137     & 0.0078     & 0.0030     \\
$\hat{D}_{CC}$              & 0.0252     & 0.0092     & 0.0063     & 0.0027   \\
$\hat{KL}$                    & 1.1060      & 0.2041     & 0.1122     & 0.0658     \\ \hline
\multicolumn{5}{c}{ {($\lambda_1$,$\lambda_2$)=(2,1)}}                                \\
                                                      \hline
Measure / $(n_1,n_2)$ & (10,10)    & (50,40)    & (80,100)   & (200,200)  \\ \hline
$\hat{D}_{Ker}$                 & 0.0209     & 0.0038     & 0.0016     & 0.0008     \\
$\hat{D}_{CC}$             & 0.0271     & 0.0096     & 0.0055     & 0.0028     \\
$\hat{KL}$                    & 0.2641     & 0.0351     & 0.0147     & 0.0077     \\ \hline
\multicolumn{5}{c}{ {($\lambda_1$,$\lambda_2$)=(1,5)}}                                \\
                                                     \hline
Measure / $(n_1,n_2)$ & (10,10)    & (50,40)    & (80,100)   & (200,200)  \\ \hline
$\hat{D}_{Ker}$                & 0.0337     & 0.0066     & 0.0038     & 0.0017     \\
$\hat{D}_{CC}$             & 0.0176     & 0.0042     & 0.0022     & 0.0012     \\
$\hat{KL}$                    & 1.0860      & 0.5165     & 0.3337     & 0.2786     \\ \hline
\multicolumn{5}{c}{ {($\lambda_1$,$\lambda_2$)=(10,9)}}                               \\
                                                      \hline
Measure / $(n_1,n_2)$ & (10,10)    & (50,40)    & (80,100)   & (200,200)  \\ \hline
$\hat{D}_{Ker}$                & 1.4509e-04 & 9.3477e-06 & 2.6519e-06 & 6.1829e-07 \\
$\hat{D}_{CC}$              & 0.0783     & 0.0209     & 0.0104     & 0.0041     \\
$\hat{KL}$                    & 0.5768     & 0.0228     & 0.0058     & 0.0019     \\ \hline
\end{tabular}
\end{table}
\FloatBarrier


\FloatBarrier
\begin{table}[]
\caption{Comparison of MSEs of Kernel-based, Empirical-based, and U-statistics-based estimators for various sample of sizes $(n_1,n_2)$ drawn from two exponential distributions with parameters ($\lambda_1$,$\lambda_2$).}
\begin{tabular}{cccccccc}
\hline
\multicolumn{8}{c}{ {($\lambda_1$,$\lambda_2$)=(0.5,0.1)}}  \\
  \hline
Estimators / $(n_1,n_2)$ & (5,10)  & (10,10) & (20,10) & (20,20) & (40,30) & (30,40) & (50,50) \\ \hline
Kernel                   & 3.380    & 3.750    & 3.440    & 1.587   & 0.977   & 0.756   & 0.621   \\
Empirical                & 4.787   & 4.133   & 3.383   & 1.700    & 0.997   & 0.789   & 0.599   \\
U-stat                   & 3.390    & 3.239   & 2.777   & 1.440    & 0.919   & 0.703   & 0.549   \\ \hline
\multicolumn{8}{c}{{($\lambda_1$,$\lambda_2$)=(0.1,0.5)} }  \\
 \hline
Estimators / $(n_1,n_2)$ & (5,10)  & (10,10) & (20,10) & (20,20) & (40,30) & (30,40) & (50,50) \\ \hline
Kernel                   & 7.565   & 3.052   & 1.592   & 1.923   & 0.804   & 1.018   & 0.644   \\
Empirical                & 6.842   & 3.328   & 1.865   & 2.031   & 0.816   & 1.011   & 0.619   \\
U-stat                   & 5.811   & 2.601   & 1.465   & 1.696   & 0.727   & 0.909   & 0.569   \\ \hline
\multicolumn{8}{c}{{($\lambda_1$,$\lambda_2$)=(0.1,1)} }   \\
  \hline
Estimators / $(n_1,n_2)$ & (5,10)  & (10,10) & (20,10) & (20,20) & (40,30) & (30,40) & (50,50) \\ \hline
Kernel                   & 8.201   & 3.527   & 1.521   & 1.567   & 0.843   & 0.996   & 0.677   \\
Empirical                & 7.576   & 3.748   & 1.671   & 1.692   & 0.873   & 1.045   & 0.694   \\
U-stat                   & 6.823   & 3.186   & 1.439   & 1.541   & 0.799   & 0.996   & 0.658   \\ \hline
\multicolumn{8}{c}{{($\lambda_1$,$\lambda_2$)=(1,0.5)} }     \\
 \hline
Estimators / $(n_1,n_2)$ & (5,10)  & (10,10) & (20,10) & (20,20) & (40,30) & (30,40) & (50,50) \\ \hline
Kernel                   & 0.1174  & 0.0758  & 0.0600    & 0.0423  & 0.0206  & 0.0199  & 0.0128  \\
Empirical                & 0.1927  & 0.0967  & 0.0639  & 0.0445  & 0.0191  & 0.0197  & 0.0122  \\
U-stat                   & 0.0936  & 0.0580   & 0.0448  & 0.0336  & 0.0162  & 0.0167  & 0.0107  \\ \hline
\multicolumn{8}{c}{{($\lambda_1$,$\lambda_2$)=(2,5)}}    \\
 \hline
Estimators / $(n_1,n_2)$ & (5,10)  & (10,10) & (20,10) & (20,20) & (40,30) & (30,40) & (50,50) \\ \hline
Kernel                   & 0.0131  & 0.0062  & 0.0033  & 0.0028  & 0.0015  & 0.0019 & 0.0011  \\
Empirical                & 0.0131  & 0.0074  & 0.0041  & 0.0029  & 0.0015  & 0.0018  & 0.0010  \\
U-stat                   & 0.0093  & 0.0048  & 0.0027  & 0.0023  & 0.0012  & 0.0015  & 0.0009 \\ \hline
\multicolumn{8}{c}{{($\lambda_1$,$\lambda_2$)=(10,5)}}   \\
\hline
Estimators / $(n_1,n_2)$ & (5,10)  & (10,10) & (20,10) & (20,20) & (40,30) & (30,40) & (50,50) \\ \hline
Kernel                   & 0.0011  & 0.0010   & 0.0006 & 0.0003 & 0.0001 & 0.0001 & 0.0001 \\
Empirical                & 0.0018 & 0.0012  & 0.0006 & 0.0004  & 0.0001 & 0.0001 & 0.0001 \\
U-stat                   & 0.0008 & 0.0007 & 0.0004 & 0.0003  & 0.0001 & 0.0001 & 0.0001 \\ \hline
\end{tabular}
\end{table}

\FloatBarrier

 1. The U-statistics-based estimator consistently achieves the lowest MSE values across different sample sizes, indicating its superior accuracy and efficiency compared to the other two estimators.
 
    2. When comparing the Kernel-based and Empirical-based estimators, their relative performance varies with sample sizes. For smaller sample sizes, the Kernel-based estimator demonstrates better accuracy, likely due to its smoothing properties, which help mitigate estimation variability. However, as the sample size increases, the Empirical-based estimator gradually outperforms the Kernel-based estimator, suggesting that it performs better for large datasets.

\subsection{Comparison of confidence intervals for $D$}

In this section, we conduct simulations to assess the performance of different confidence intervals for the measure \( D \), as discussed in Sections 2.4–2.5. Each simulation is repeated 1,000 times. To compare the effectiveness of these confidence intervals, we consider the following key criteria:

\begin{itemize}
    \item Coverage Probability (CP): This measures the proportion of times the true parameter value falls within the constructed confidence interval (CI). Ideally, the empirical coverage probability should be close to the nominal confidence level. A lower coverage error—defined as the absolute difference between the observed coverage probability and the nominal level—indicates a more reliable method.
\vspace{2mm}

    \item Average Length (AL) of CIs: A shorter confidence interval suggests higher precision in estimating the unknown parameter. Therefore, methods yielding confidence intervals with smaller average lengths are generally preferred, provided they maintain adequate coverage probability.  
\end{itemize}

\vspace{2mm}


 For various sample sizes and parameter values, we calculated the CPs and ALs of confidence intervals obtained based on the Normal approximation and JEL methods (given by \eqref{eq:Normal App} and\eqref{eq: CI JEL}, respectively)  . These results are presented in Table 3. We implement the jackknife empirical likelihood (JEL)-based confidence intervals using the R package \texttt{empilik}. The following observations can be drawn from Table 3:

1. For various configurations of the exponential model parameters, the average lengths (ALs) of confidence intervals based on the Jackknife Empirical Likelihood (JEL) and Normal approximation methods decrease as the sample size increases. Additionally, the coverage probabilities (CPs) of both methods converge toward the nominal level of 95\% as the sample size increases.

2. Table 3 shows that the JEL method consistently achieves coverage probabilities closer to 95\% across different conditions compared to the Normal approximation. Moreover, the JEL method produces significantly shorter confidence intervals than the Normal approximation, highlighting its efficiency.

From this study, we conclude that the JEL method outperforms the Normal approximation for constructing confidence intervals for the divergence measure $D$. Furthermore, Theorem 3 indicates that the JEL method can be effectively applied to hypothesis testing for assessing the significance of the divergence measure.

\FloatBarrier
	\begin{table}[]
		\caption{Coverage Probability (CP) and Average length (AL) of confidence intervals for measure $D$ at 95\% confidence level}
		\begin{tabular}{ll  ll  ll  }
			\hline
			\multicolumn{1}{c}{\multirow{2}{*}{$(\lambda_1,\lambda_2)$}} & \multicolumn{1}{c}{\multirow{2}{*}{$(n_1,n_2)\;\;\;$}} & \multicolumn{2}{c}{JEL}  & \multicolumn{2}{c}{Normal}  \\ \cline{3-4} \cline{5-6} 
			\multicolumn{1}{c}{}                                     & \multicolumn{1}{c}{}                            & \multicolumn{1}{c}{CP (\%)} & \multicolumn{1}{c}{AL} & \multicolumn{1}{c}{CP (\%)} & \multicolumn{1}{c}{AL} \\ \hline
			
			(0.1,0.1)	& (10,10)                                                &      92                     &    5.215   & 94 & 4.154        \\
			&    (30,40)                                           &     93                        &     1.372     & 93 & 1.547           \\
			&  (70,50)                                               &  96                          &    0.981        &96 & 0.975      
            \\
            
			&  (100,100)                                               &     96                       &     0.952    & 96 & 0.961         
            
            \\ \hline
			(0.1,0.5)	& (10,10)                                                &    82                        &       5.140   & 81 & 6.120               \\
			&      (30,40)                                           &    87                        &     4.336        & 84 & 4.197          \\
			&  (70,50)                                               &     91                        &     2.787        & 85 & 2.882             
            \\
            
			&  (100,100)                                               &     92                       &     2.162    & 89 & 2.287         
            \\\hline
			(1,0.5)	& (10,10)                                                &    77                         &       1.059   & 79 & 0.986               \\
			&      (30,40)                                           &    86                       &     0.554         & 86 & 0.541         \\
			&  (70,50)                                               &     92                       &     0.316        & 91 & 0.301             
            \\
            
			&  (100,100)                                               &     94                       &     0.262    & 92 & 0.267     
			\\\hline
			(1,2)	& (10,10)                                                &    80                         &       0.780           & 85 & 0.460      \\
			&      (30,40)                                           &    84                     &     0.315  & 85 & 0.397                  \\
			&  (70,50)                                               &     88                      &     0.206         &     89                       &     0.212      \\
            
			&  (100,100)                                               &     91                     &     0.172    & 91 & 0.181         
			\\\hline
			(10,1)	& (10,10)                                                &                          84   &          0.724 & 84 & 0.702             \\
			&      (30,40)                                           &    95                       &       0.421   & 93 & 0.366             \\
			&  (70,50)                                               &     96                     &     0.215             & 96 & 0.211   \\
            
			&  (100,100)                                               &     95                     &     0.112    & 95 & 0.124            \\\hline
		\end{tabular}
		
	\end{table}

\FloatBarrier

\section{Real Data Analysis}

In this section, we apply the divergence measure $D$ to real-world datasets to demonstrate its practical utility. We focus on datasets containing uncensored cases where survival analysis is relevant, and we assess the ability of this new measure to distinguish between two populations. To visually examine differences in survival functions, we plot the empirical survival curves for the groups under comparison. These plots allow us to qualitatively assess how the survival probabilities of the two groups evolve over time. Furthermore, we evaluate whether the values obtained from the divergence measure correspond to the observed differences between the curves. Through these analyses, we illustrate the effectiveness of our measure in capturing and quantifying differences between survival functions.

\noindent \textbf{Veteran Dataset:} The Veteran dataset is a publicly available dataset on \href{https://www.kaggle.com/code/aradinka/survival-analysis-veterans-lung-cancer-study}{Kaggle}\footnote{\url{https://www.kaggle.com/code/aradinka/survival-analysis-veterans-lung-cancer-study}}. This dataset is a well-known survival analysis dataset derived from a clinical trial on patients with advanced lung cancer. It contains variables such as treatment type (standard vs. test treatment), survival time, and patient characteristics, including age, cell type, and performance status.

We restricted this dataset to uncensored cases in order to examine differences in survival functions between the two treatment groups. In this analysis, the value of the divergence measure $D$ is \textbf{2.3972}. Since this value is relatively low, it suggests that there is not a substantial difference in the survival functions of patients receiving the standard treatment compared to those receiving the test treatment. Notably, the result of the measure $D$ is consistent with the visual comparison of the empirical survival curves, as shown in Figure 1.

\FloatBarrier
\begin{figure}[h]
  \centering
  \includegraphics[width=1\textwidth]{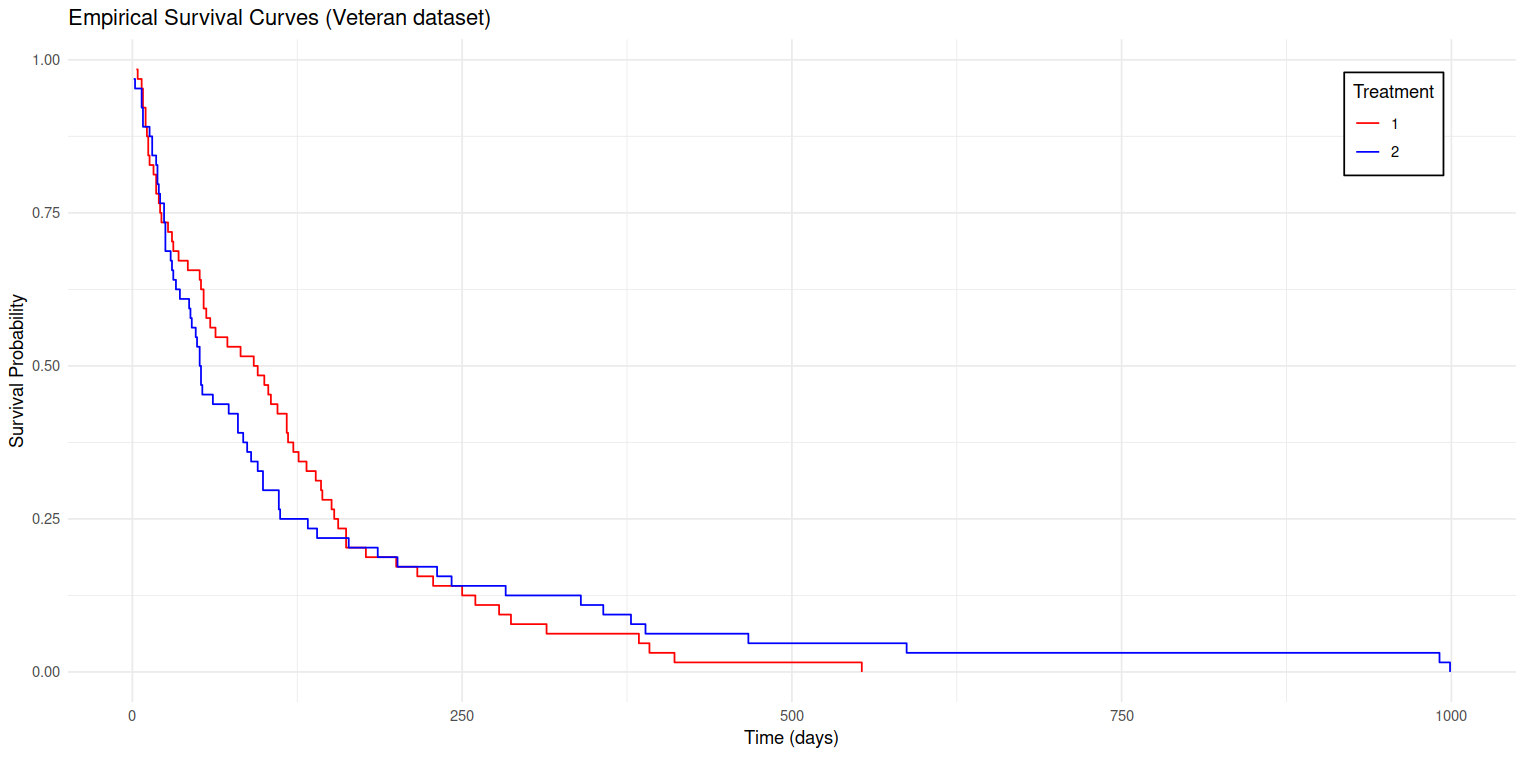}
  \caption{Empirical survival curves for comparing the survival functions of the standard and test treatment groups in the Veteran dataset.}
  \label{fig:example}
\end{figure}

\noindent \textbf{Lung Cancer Dataset:}  The Lung Cancer dataset is a publicly available dataset on \href{https://www.kaggle.com/datasets/mysarahmadbhat/lung-cancer}{Kaggle}\footnote{\url{https://www.kaggle.com/datasets/mysarahmadbhat/lung-cancer}}. This dataset is a widely used dataset in survival analysis, comprising information on patients diagnosed with advanced lung cancer. It includes survival time, censoring status, and relevant covariates such as sex, age, and treatment type (standard vs. test treatment).

In this analysis, we focus on uncensored cases and assess the differences in survival functions between the two sex groups. The value of the divergence measure is \textbf{6.9714}, which is high and suggests a substantial difference in survival functions between male and female patients. This observation is consistent with the separation evident in the empirical survival curves, as shown in Figure 2.\\
\FloatBarrier
\begin{figure}[h]
  \centering
  \includegraphics[width=1\textwidth]{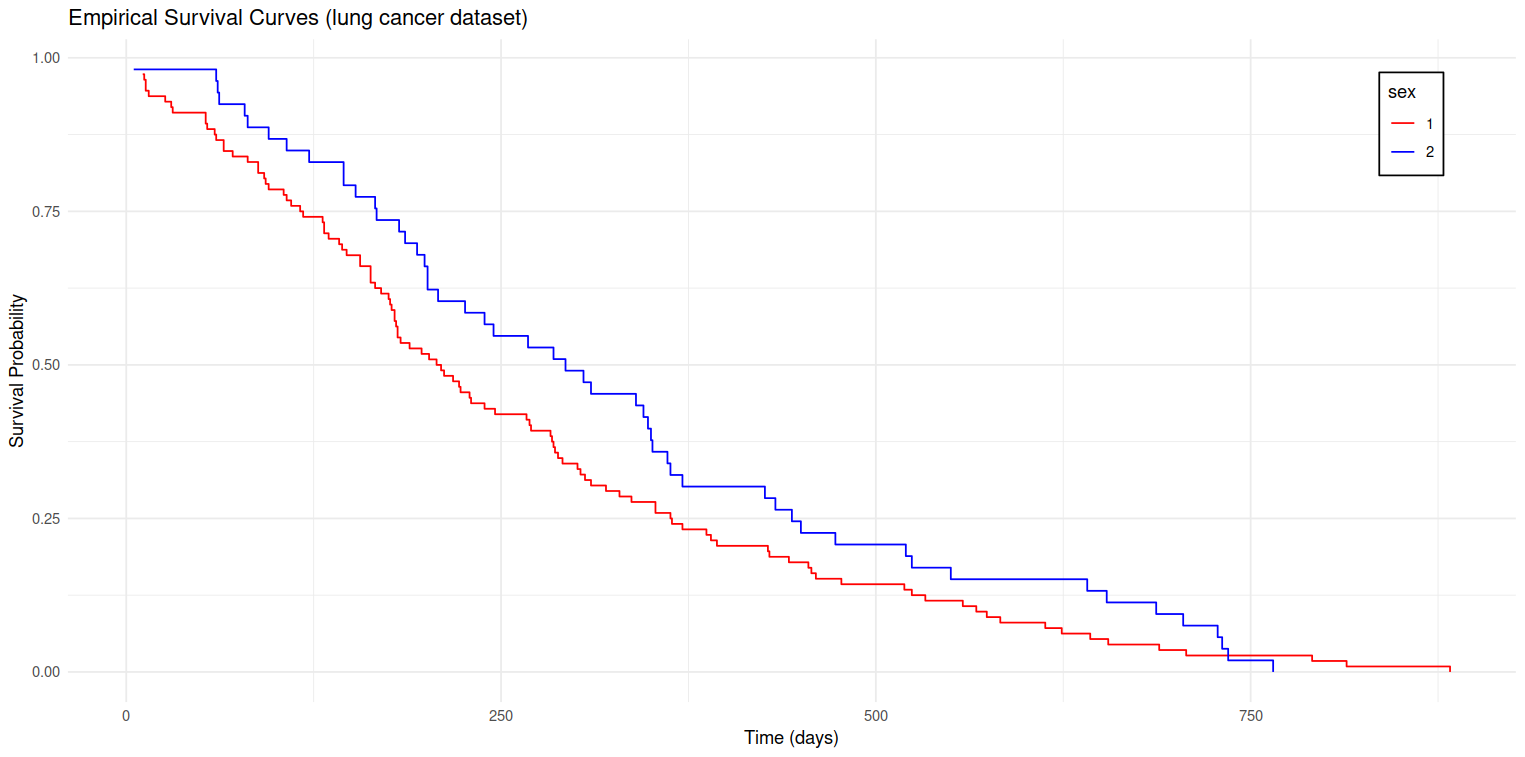}
  \caption{Empirical survival curves for comparing the survival functions of male and female patients in the Lung Cancer dataset.}
  \label{fig:example}
\end{figure}
\FloatBarrier

\noindent \textbf{PBC Dataset:} The PBC dataset is a publicly available dataset on \href{https://www.kaggle.com/datasets/homayoonkhadivi/primary-biliary-cirrhosis-pbc-disease-dataset}{Kaggle}\footnote{\url{https://www.kaggle.com/datasets/homayoonkhadivi/primary-biliary-cirrhosis-pbc-disease-dataset}}. This dataset is a well-known dataset in survival analysis, originating from a clinical trial conducted by the Mayo Clinic on patients diagnosed with primary biliary cirrhosis, a chronic liver disease. The dataset includes variables such as survival time, censoring status, sex, treatment assignment, and clinical covariates like age, bilirubin levels, albumin, and prothrombin time.

In this analysis, we assess the differences in survival functions between the two sex groups, considering only uncensored cases. The divergence measure is \textbf{22.9262}, indicating a substantial difference in survival outcomes between male and female patients. This finding aligns with the separation observed in the empirical survival curves, as illustrated in Figure 3.

\FloatBarrier
\begin{figure}[h]
  \centering
  \includegraphics[width=1\textwidth]{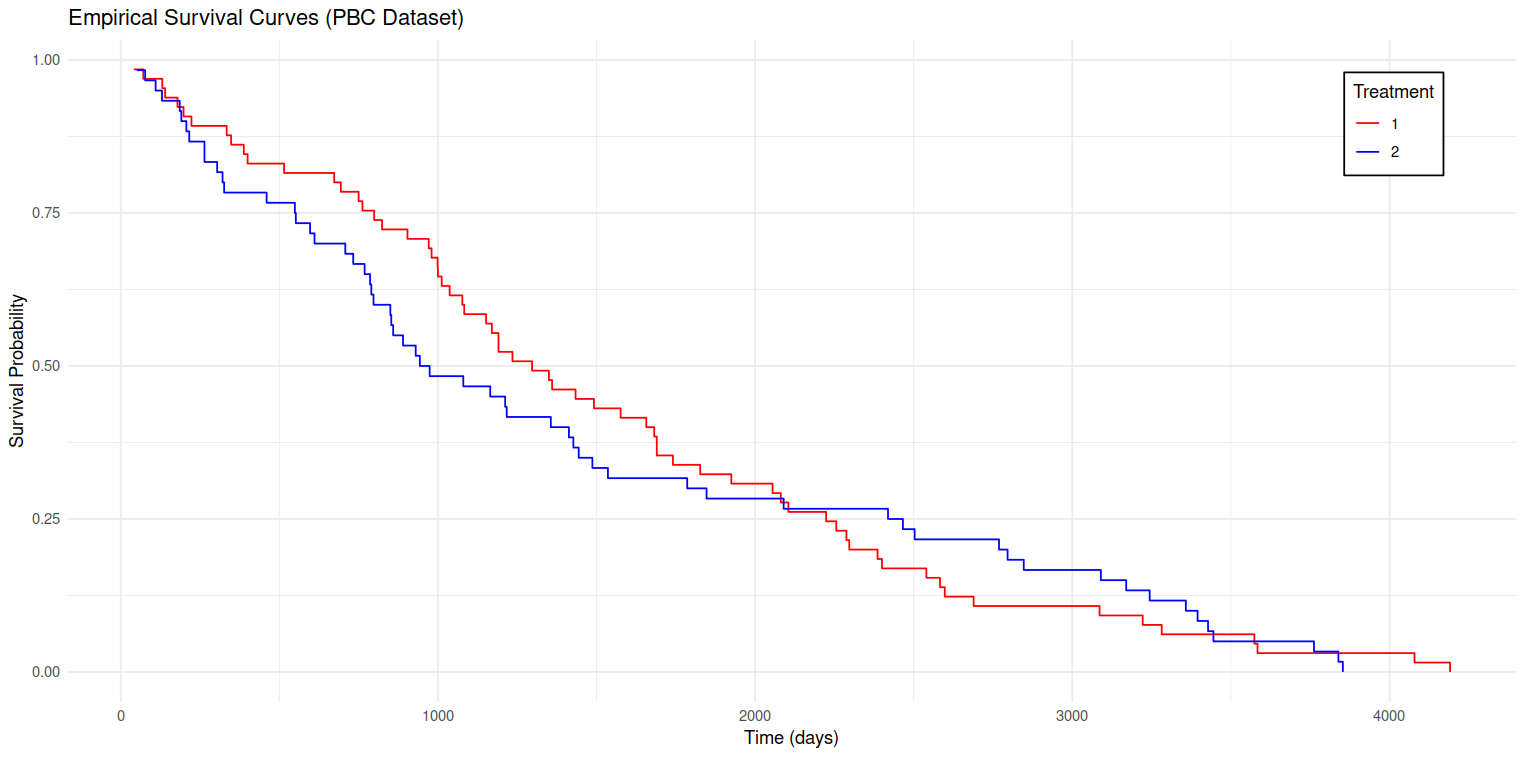}
  \caption{Empirical survival curves for comparing the survival functions of male and female patients in the PBC dataset.}
  \label{fig:example}
\end{figure}
\FloatBarrier

\noindent \textbf{GBSG2 dataset:} The GBSG2 dataset is a publicly available dataset on \href{https://r-packages.io/datasets/GBSG2}{R}\footnote{\url{https://r-packages.io/datasets/GBSG2}}. This dataset originates from a study conducted by the German Breast Cancer Study Group and contains information on women diagnosed with breast cancer. The dataset includes survival time, censoring status, and clinical variables such as age, tumor size, number of positive lymph nodes, hormone receptor status, and treatment information, including hormonal therapy.

In this analysis, we focus exclusively on uncensored cases to evaluate differences in survival functions between patients receiving hormonal therapy and those who did not. The divergence measure is \textbf{7.5641}, reflecting a notable difference in survival outcomes between the two groups. This result is supported by the separation seen in the empirical survival curves, as depicted in Figure 4.

\FloatBarrier
\begin{figure}[h]
  \centering
  \includegraphics[width=1\textwidth]{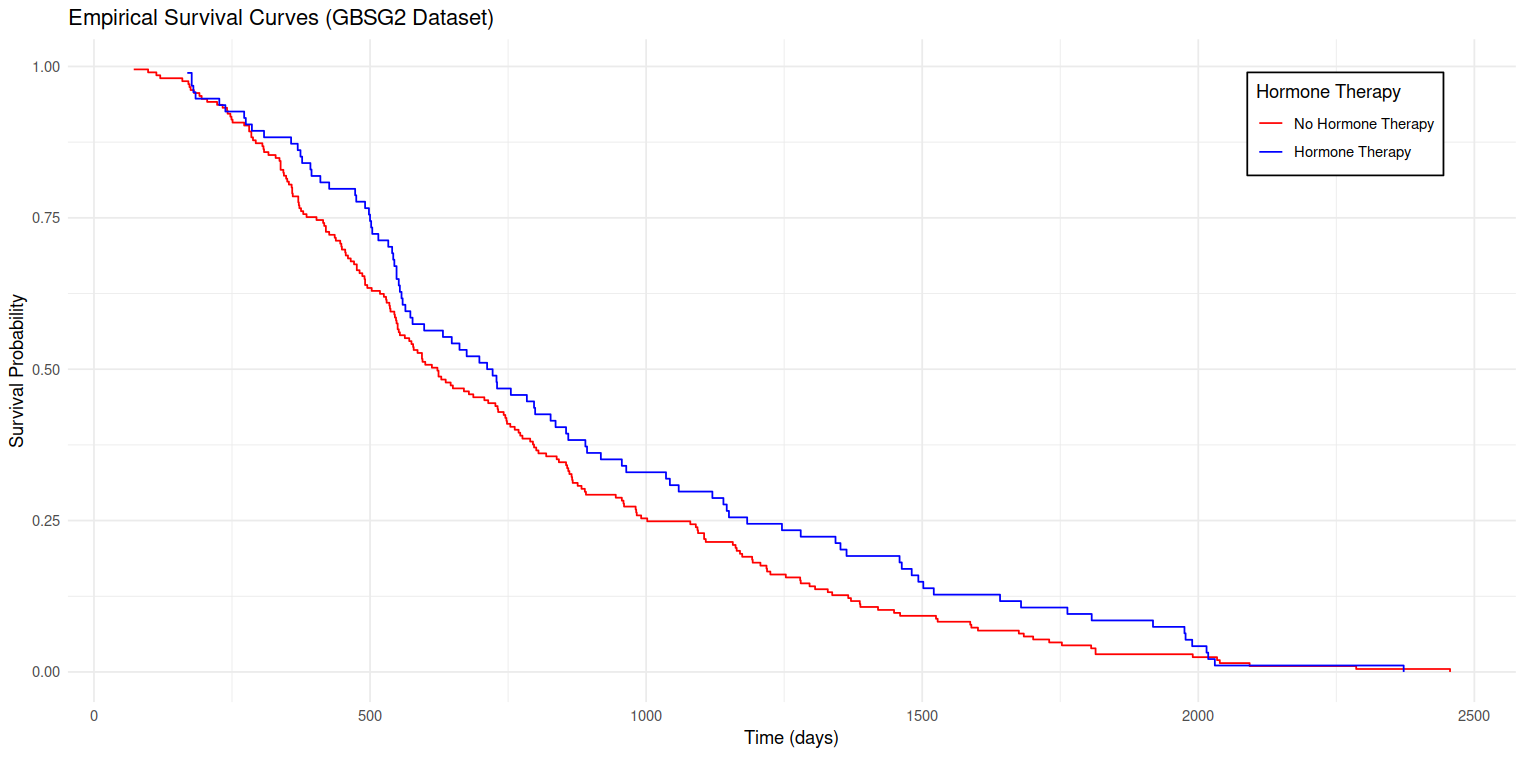}
  \caption{Empirical survival curves for comparing the survival functions of patients receiving hormonal therapy versus no hormonal therapy in the GBSG2 dataset.}
  \label{fig:example}
\end{figure}
\FloatBarrier

\subsection{Application on Image Datasets}

In this section, we use image datasets to further demonstrate the applications of the divergence measure. To ensure consistency across the data, the selected images were resized to a uniform dimension and normalized to standardize pixel intensities.

\noindent \textbf{Brain Tumor MRI dataset:}

The Brain Tumor MRI Images dataset is a publicly available dataset on \href{https://www.kaggle.com/datasets/sartajbhuvaji/brain-tumor-classification-mri/data}{Kaggle}\footnote{\url{https://www.kaggle.com/datasets/sartajbhuvaji/brain-tumor-classification-mri/data}}. It contains MRI images of human brains divided into four categories: No Tumor, Pituitary Tumor, Benign Tumor, and Malignant Tumor. For our study, we selected three categories: No Tumor (NT), Benign Tumor (BT), and Malignant Tumor (MT). From each category, we chose three MRI images representing different parts of the brain that best illustrate the characteristics of each class. The categories included in this study are described as follows:

The No Tumor (NT) category contains MRI images of healthy brains with no abnormal growths or masses. In these images, the brain tissue appears normal, and there are no signs of tumors or unusual structures. These images serve as a baseline, allowing for effective comparison and detection of abnormalities in other MRI scans.

\begin{figure}[h]
  \centering

  \begin{subfigure}[t]{0.32\textwidth}
    \centering\includegraphics[width=\linewidth,height=4cm,keepaspectratio]{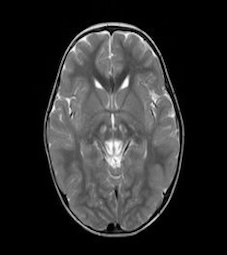}
    \caption{NT1}
  \end{subfigure}%
  \begin{subfigure}[t]{0.32\textwidth}
    \centering\includegraphics[width=\linewidth,height=4cm,keepaspectratio]{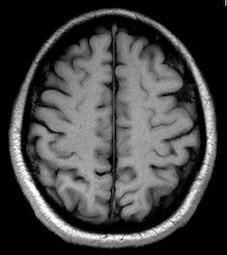}
    \caption{NT2}
  \end{subfigure}%
  \begin{subfigure}[t]{0.32\textwidth}
    \centering\includegraphics[width=\linewidth,height=4cm,keepaspectratio]{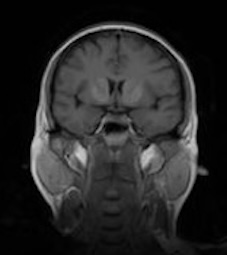}
    \caption{NT3}
  \end{subfigure}

  \vspace{2pt} 

  \begin{subfigure}[t]{0.32\textwidth}
    \centering\includegraphics[width=\linewidth,height=4cm,keepaspectratio]{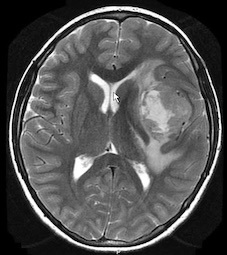}
    \caption{BT1}
  \end{subfigure}%
  \begin{subfigure}[t]{0.32\textwidth}
    \centering\includegraphics[width=\linewidth,height=4cm,keepaspectratio]{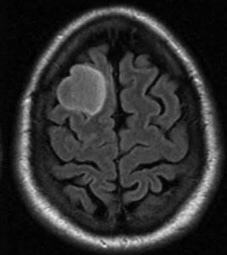}
    \caption{BT2}
  \end{subfigure}%
  \begin{subfigure}[t]{0.32\textwidth}
    \centering\includegraphics[width=\linewidth,height=4cm,keepaspectratio]{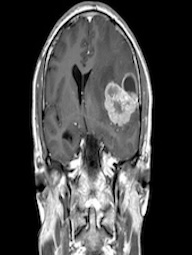}
    \caption{BT3}
  \end{subfigure}

  \vspace{2pt}

  \begin{subfigure}[t]{0.32\textwidth}
    \centering\includegraphics[width=\linewidth,height=4cm,keepaspectratio]{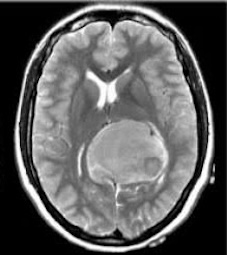}
    \caption{MT1}
  \end{subfigure}%
  \begin{subfigure}[t]{0.32\textwidth}
    \centering\includegraphics[width=\linewidth,height=4cm,keepaspectratio]{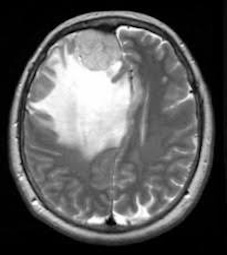}
    \caption{MT2}
  \end{subfigure}%
  \begin{subfigure}[t]{0.32\textwidth}
    \centering\includegraphics[width=\linewidth,height=4cm,keepaspectratio]{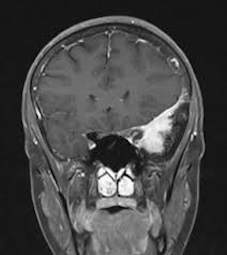}
    \caption{MT3}
  \end{subfigure}

  \caption{Axial and Coronal MRI Images of Brains Representing Different Tumor Categories}
  \label{fig:mri_grid_3x3}
\end{figure}

The Benign Tumor (BT) category includes MRI images showing non-cancerous tumors in the brain. Benign tumors grow slowly and do not spread to other parts of the brain or body. While generally less dangerous than malignant tumors, they can still cause symptoms if they press on surrounding brain tissue. Accurate identification of benign tumors is important for planning appropriate treatment, which may involve monitoring, minor surgery, or other non-aggressive interventions.
\\

The Malignant Tumor (MT) category consists of MRI images showing cancerous tumors. These tumors are aggressive, can grow rapidly, and may spread to other parts of the brain or body. Malignant tumors require immediate and intensive treatment, such as surgery, chemotherapy, or radiotherapy. Differentiating malignant tumors from benign tumors is crucial, as it directly impacts the treatment strategy and can significantly affect patient outcomes.

Differentiating these categories is essential for early detection, accurate diagnosis, and effective treatment planning. By distinguishing between healthy brains, benign tumors, and malignant tumors, medical professionals can identify abnormalities more reliably, reduce misdiagnoses, and streamline medical care. In this study, we apply our divergence measure to differentiate between the selected MRI images, as illustrated in Figure 5. In Figure 5, the first row shows MRI images from NT patients, the second row shows images from BT patients, and the third row shows images from MT patients.

Using the divergence measure, we calculated estimates for each pair of MRI images and compared the three categories. For calculation, each grayscale image is represented by its pixel intensity values (scaled to [0,1]). These values are flattened into a one dimensional vector and used as the sample data for subsequent analysis. For any two images, denoted data1 and data2, we then computed an empirical based divergence estimate.
The resulting values are presented in Tables 4, 5, and 6. A higher value of the divergence estimate indicates a greater likelihood of a tumor. These results demonstrate that the divergence measure can serve as an effective tool for detecting tumors in medical image data.

While our study demonstrates the effectiveness of the divergence measure in differentiating brain MRI images across tumor categories, the results could be further strengthened with longitudinal data, where MRI scans of the same individuals are available at different stages of tumor progression, captured with consistent size and image quality. The current dataset, being cross-sectional, limits the ability to draw more definitive conclusions regarding tumor evolution. Incorporating longitudinal data in future studies would provide deeper insights and allow for more robust validation of the divergence measure in tracking disease progression.

\FloatBarrier
\begin{table}[h]
\caption{Estimated Divergence Measure Values for Pairs of Different Classes in Axial View Brain MRIs}
\begin{tabular}{l|lll}
    & NT1   & BT1   & MT1   \\ \hline
NT1 & 0     & 0.037 & 0.049 \\
BT1 & 0.037 & 0     & 0.012 \\
MT1 & 0.049 & 0.021 & 0    
\end{tabular}
\end{table}

\begin{table}[h]
\caption{Estimated Divergence Measure Values for Pairs of Different Classes in Axial View Brain MRIs}
\begin{tabular}{l|lll}
    & NT2   & BT2   & MT2   \\ \hline
NT2 & 0     & 0.024 & 0.093 \\
B2 & 0.024 & 0     & 0.046 \\
MT2 & 0.093 & 0.046 & 0    
\end{tabular}
\end{table}

\begin{table}[h]
\caption{Estimated Divergence Measure Values for Pairs of Different Classes in Coronal View Brain MRIs}
\begin{tabular}{l|lll}
    & NT2   & BT2   & MT2   \\ \hline
NT2 & 0     & 0.023 & 0.032 \\
B2 & 0.023 & 0     & 0.023 \\
MT2 & 0.032 & 0.023 & 0    
\end{tabular}
\end{table}
\FloatBarrier

\noindent \textbf{Vehicle Detection Image Dataset:}\\
The Vehicle Detection Image Dataset is a publicly available dataset on \href{https://www.kaggle.com/datasets/pkdarabi/vehicle-detection-image-dataset?utm_source=chatgpt.com}{Kaggle}\footnote{\url{https://www.kaggle.com/datasets/pkdarabi/vehicle-detection-image-dataset?utm_source=chatgpt.com}}. It contains high-resolution images captured from various urban areas, showing different types of vehicles such as cars, trucks, buses, and motorcycles on roads. For our study, we selected three similar-looking images of a road taken at different times, as shown in Figure 6. For calculation, each colored image was converted to grayscale by averaging the red, green, and blue channel intensities at each pixel. The resulting grayscale image was then flattened into a one dimensional vector, which was used as the sample data for subsequent analysis.

The results show that, even though the images look almost the same at first glance, the divergence measure can detect small differences between them. This example further shows that the divergence measure can be a useful tool for finding small changes in colored images as well.

\begin{figure}[h]
  \centering

  \begin{subfigure}[t]{0.32\textwidth}
    \centering\includegraphics[width=\linewidth,height=4cm,keepaspectratio]{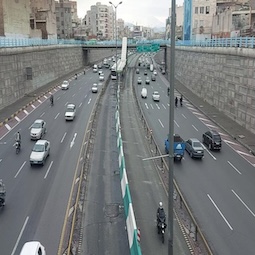}
    \caption{T1}
  \end{subfigure}%
  \begin{subfigure}[t]{0.32\textwidth}
    \centering\includegraphics[width=\linewidth,height=4cm,keepaspectratio]{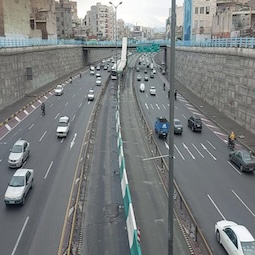}
    \caption{T2}
  \end{subfigure}%
  \begin{subfigure}[t]{0.32\textwidth}
    \centering\includegraphics[width=\linewidth,height=4cm,keepaspectratio]{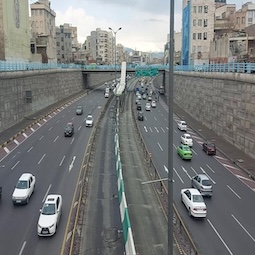}
    \caption{T3}
  \end{subfigure}

  \vspace{2pt} 

  \caption{Images of Roads Captured at Different Times Showing Vehicles}
  \label{fig:mri_grid_3x3}
\end{figure}

\begin{table}[h]
\caption{Estimated Divergence Measure Values for Pairs of Road Images Taken at Different Times}
\begin{tabular}{l|lll}
    & T1   & T2   & T3   \\ \hline
T1 & 0     & 0.000037 & 0.004760 \\
T2 & 0.000037 & 0     & 0.004980 \\
T3 & 0.004760 & 0.004980 & 0    
\end{tabular}
\end{table}

\section{Concluding Remarks}

In this paper, we considered a divergence measure $D$ based on the squared difference between two survival functions. We provided several motivations for this measure and established its connection to extropy. Three non-parametric estimators of $D$ were derived, and through a simulation study, we observed that the U-statistic-based estimator outperforms the other two.\\

We also compared the measure $D$ with existing divergence measures, such as \( D_{CC} \) (introduced by \cite{cox2016practical}) and the Kullback–Leibler (KL) divergence. Simulation results indicate that the estimator of $D$ often outperforms the estimators of existing measures in various scenarios under the mean squared error (MSE) criterion. Moreover, desirable properties such as symmetry, positivity, differentiability, and smoothness further support the utility of this measure $D$.

Additionally, we addressed the problem of constructing confidence intervals for the measure $D$. We employed both the normal approximation for U-statistics and the jackknife empirical likelihood (JEL) approach. Our simulation study showed that the JEL method outperforms the normal approximation, particularly when the data are generated from two exponential distributions.

Overall, our findings suggest that the divergence measure $D$ is a promising alternative for analyzing data, with robust estimation methods and reliable inference procedures. As part of our future work, we are developing non-parametric estimators of $D$ when survival times are right censored. 

We have also proposed a general divergence measure $D_w$ (defined in \eqref{eq:1.3}). There are several advantages in using \( D_w \).  Its smoothness and differentiability facilitate statistical inference and optimization. It is designed to handle lifetime data and right-censored data effectively.  The link to extropy provides a robust theoretical basis.  Finally, squared integral measures like \( D_w \) often exhibit desirable statistical properties, such as reduced bias and improved convergence. In future work, we aim to develop non-parametric estimators and confidence intervals for this generalized measure \( D_w \) as well for various choices of the weight function $w$.


\bibliographystyle{apalike}
\bibliography{references}		


\end{document}